\documentclass[12pt]{article}
\usepackage{amsfonts}
\usepackage{amssymb}
\usepackage{amscd}
\usepackage{amsmath}
\usepackage{epsfig}
\usepackage{graphicx, subfigure}
\usepackage{latexsym}
\usepackage{verbatim}
\usepackage{float}
\usepackage{hyperref}
\usepackage{multirow}

\newlength{\dinwidth}
\newlength{\dinmargin}
\begin{document}

\title{From formal to actual Puiseux series solutions of algebraic differential equations
of first order}

\author{Vladimir Dragovi\'c$^1$, Renat Gontsov$^2$, and Irina Goryuchkina$^3$}

\date{}

\maketitle

\footnotetext[1]{Department of Mathematical Sciences, University
	of Texas at Dallas, 800 West Campbell Road, Richardson TX 75080,
	USA. Mathematical Institute SANU, Kneza Mihaila 36, 11000
	Belgrade, Serbia.  E-mail: {\tt Vladimir.Dragovic@utdallas.edu}}

\footnotetext[2]{M.S. Pinsker Laboratory no.1, Institute for Information Transmission Problems of the Russian Academy of Sciences, Bolshoy Karetny per. 19, build.1, Moscow 127051 Russia.   E-mail: {\tt gontsovrr@gmail.com}}

\footnotetext[3]{Keldysh Institute of Applied Mathematics, Russian Academy of Sciences, Miusskaya sq. 4, Moscow 125047
Russia. E-mail: {\tt  igoryuchkina@gmail.com}}

\begin{abstract}
The existence, uniqueness and convergence of formal Puiseux series solutions of non-autonomous algebraic differential equations of first order at a nonsingular point of the equation is studied, including the case where the celebrated Painlev\'e theorem cannot be applied explicitly for the study of convergence. Several examples illustrating relationships to the Painlev\'e theorem and lesser-known Petrovi\'c's results are provided.
\end{abstract}


\section{Introduction}

The solution of the algebraic equation
\begin{equation}\label{alg}
F(x,y)=\sum\limits_{i=1}^na_i\,x^{p_i}y^{q_i}=0, \qquad x, y\in\mathbb C,
\end{equation}
is a unique object, an algebraic function
 $y=y(x)$ if the polynomial $F$ is irreducible or a set of algebraic functions if $F$ is reducible. When one talks about local properties of the solution, then in a  neighborhood of almost every point $x=x_0\in{\mathbb C}$ there are finitely many holomorphic germs of the function $y(x)$, whereas in a neighborhood of each point $x=x_0$ from the finite set of ramification points there are finitely many germs of this function which can be presented as convergent Puiseux series in fractional powers of the variable  $x-x_0$. The construction of the Newton--Puiseux polygon of an algebraic equation allows to find all the germs. One more remark is that every Puiseux series, which formally satisfies the equation (\ref{alg}), has a nonzero radius of convergence, see \cite{Ploski} for a contemporary account and \cite[Ch. 4.3]{Walker} for a more classical approach.
\smallskip

When we consider an algebraic {\it differential} equation
\begin{equation}\label{eq1}
F(x,y,y')=\sum\limits_{i=1}^na_i(x)\,y^{p_i}(y')^{q_i}=0,
\end{equation}
$a_i$ being polynomials, then the situation is quite different. The solution is not a unique object any more. In such a case, we consider an {\it entire collection of solutions} gathered under one notion -- the {\it general solution} of the equation. In general, such solutions may have singularities of non-algebraic type, and also formal power series in the variable  $x-x_0$, which satisfy (\ref{eq1}), may diverge. As a simple example, we consider the equation
$$
xy'-1=0.
$$
Its general solution possesses a singularity of non-algebraic type at the point $x=0$. Another, notable example is the Euler equation (see \cite[Ch. II]{Ha}),
$$
x^2y'-y+x=0,
$$
which has a formal solution in a form of the power series
$$
\sum_{k=0}^{\infty}k!x^{k+1}
$$
in the variable $x$ with the radius of convergence equal to zero.
\medskip

We come to two natural and important questions:
\begin{itemize}
\item [i)] {\it For which points $x=x_0\in\mathbb C$ one can guarantee the algebraic local behaviour of the general solution of the equation $(\ref{eq1})$?}
\item [ii)] {\it For which Puiseux series in the variable $x-x_0$ formally satisfying the equation $(\ref{eq1})$, one can guarantee convergence?}
\end{itemize}
(Throughout the paper we always mean by $x_0$ a finite point of $\mathbb C$. We don't speak about the infinite point separately, since it is studied, as usual, {\it via} the change $t=1/x$ of the independent variable and considering the point $t=0$ of the transformed equation.) These two questions are related to each other in fact, by mean of the ''fundamental existence theorem'' \cite{RS}, \cite{Ra}. The latter states that for each formal power series solution $\varphi\in{\mathbb C}[[x-x_0]]$ of (\ref{eq1}) there exists an actual solution having $\varphi$ as an asymptotic expansion in some sector with the vertex at $x_0$. Hence the existence of a {\it divergent} Puiseux series in the variable $x-x_0$ formally satisfying the equation (\ref{eq1}) implies that the equation possesses a solution for which $x=x_0$ is a singular point of a non-algebraic type.

At the same time, the answer to Question i) partially follows from the celebrated Painlev\'e theorem \cite{Painleve}, \cite{Painleve0}, \cite{Painleve1}. According to this theorem every solution to the equation (\ref{eq1}) can only have a singularity of an {\it algebraic type} at any point $x=x_0$, {\it with the exception of points of some fixed finite set $\Sigma$ at most determined by the equation}. In other words, non-algebraic singular points of the solutions of a first order algebraic ODE cannot fill domains in $\mathbb C$. Therefore, any Puiseux series in the variable $x-x_0$ formally satisfying the equation (\ref{eq1}) converges near $x_0$, if $x_0\in{\mathbb C}\setminus\Sigma$, which partially answers Question ii). In particular, for an {\it autonomous} algebraic ODE of first order the set $\Sigma$ is empty and hence any formal Puiseux series solution of such an equation has a nonzero radius of convergence (although after the change $t=1/x$ the transformed equation becomes non-autonomous, the last statement on convergence is still true for formal Puiseux series solutions in the variable $1/x$ considered in a neighbourhood of infinity, which was proved in \cite{CFS2019}).

However, the question of detecting the set $\Sigma$ of potentially non-algebraic singular points of the general solution of the equation (\ref{eq1}) is quite delicate and elaborate. Writing the equation for a moment as
\begin{equation}\label{eq11}
F(x,y,y')=A_0(x,y)+A_1(x,y)y'+\ldots+A_s(x,y)(y')^s=0,
\end{equation}
one definitely says that $\Sigma$ contains the points $x=x_0$ for which
\begin{itemize}
	\item [-] either $A_s(x_0,y)\equiv0$;
	\item [-] or the equations $A_0(x_0,y)=0$, $A_1(x_0,y)=0,\ldots,A_s(x_0,y)=0$ have a common solution;
	\item [-] or after the change of variable $y=1/w$, for the new equation of the form (\ref{eq11}) with the coefficients $\widetilde A_i(x,w)$, one has $\widetilde A_0(x_0,0)=\widetilde A_1(x_0,0)=\ldots=\widetilde A_s(x_0,0)=0$.
\end{itemize}
Yet $\Sigma$ may contain additional points. In the case when the system
\begin{equation}\label{systP}
F(x,y,y')=0, \quad \frac{\partial F}{\partial y'}(x,y,y')=0, \quad
\frac{\partial F}{\partial x}(x,y,y')+y'\,\frac{\partial F}{\partial y}(x,y,y')=0
\end{equation}
has a finite set of isolated solutions $(x_0,y_0,y'_0)$, such additional points are those $x=x_0$ coming from these solutions, see \cite[p. 42]{Painleve} (the explanation is also given by Picard \cite[Ch. II, pp. 40--41]{Picard} based on the
results of Briot and Bouquet \cite{BB}).

But in the case when the Painlev\'e--Picard system (\ref{systP}) is compatible for every $x\in\mathbb C$, detecting the set $\Sigma$ is not clarified in general and thus Questions i), ii) have no an immediate answer. Moreover, even when one would succeed in obtaining the set $\Sigma$, for its points $x_0$ some Puiseux series in the variable $x-x_0$ formally satisfying the equation (\ref{eq1}) could also converge (see Example 2 in the last section). Therefore our aim is to present some general analysis of the question of convergence without appealing to the set $\Sigma$. We indeed propose such an analysis for any {\it nonsingular} point $x_0$ of the equation (\ref{eq1}), in the sense of the following definition.
\medskip

{\bf Definition 1.} A point  $x=x_0\in \mathbb C$ will be referred to as  a {\it singular point of the equation} (\ref{eq1}) if it is a zero of any of the coefficients $a_i$.
\medskip

We prove that every Puiseux series in powers of the variable $x-x_0$, starting with a generic term $c(x-x_0)^{\lambda}$, $\lambda\in{\mathbb Q}^*$, and formally satisfying (\ref {eq1}), converges in a neighborhood of the nonsingular point $x=x_0$ of the equation (\ref {eq1}). The last statement is content of Theorem 2.
\smallskip

Before we investigate the question of convergence of formal Puiseux series satisfying a {\it non-autonomous} algebraic differential equation of first order,  we will study the existence and uniqueness of such formal Puiseux series solutions in Theorem 1. To that end we will employ the Newton--Puiseux polygons in the form in which it was proposed and applied by Mihailo Petrovi\'c to study the local behavior of solutions in a neighborhood of a {\it nonsingular} point of the equation. Thus, the Petrovi\'c method is exposed in the next section. As a side remark let us also mention that Petrovi\'c observed an interesting application of {\it non-autonomous} algebraic differential equations of first order in chemical dynamics in \cite{Petrovich3}.

\section{Petrovi\'c's polygonal method and theorems}

The Serbian mathematician Mihailo Petrovi\'c Alas, a student of Emile Picard and Charles Hermite, defended his thesis \cite{Petrovich1} in 1894.  One of the chapters of his thesis is devoted to the study of the analytic properties of solutions of first-order algebraic differential equations. He developed a method, the {\it Petrovi\'c polygonal method}, applicable to algebraic differential equations of any order. His method uses the same principles as the Newton--Puiseux polygonal method for algebraic equations. Moreover, the Petrovi\'c polygon differs from the polygons of C. Briot and J. Bouquet, and of H.  Fine \cite{Fine}, who also generalized the Newton--Puiseux polygonal method. In \cite{DGAHES2020}, a comparative analysis of constructions of polygons of Petrovi\'c and Fine was performed and Theorem 7 therein establishes
the relationship between the two methods under certain conditions, see also \cite{DGBAMS2020}, Theorem 2, where the notion of ''the Fine--Petrovi\'c'' polygons was coined.
\smallskip

Petrovi\'c investigated the local behavior of the general solution of the equation \eqref {eq1} in the neighborhood of its {\it nonsingular} point $x = x_0$. According to the Painlev\'e theorem, the solutions of the equation \eqref {eq1} generically  have the local form $y = (x-x_0)^{\lambda}f(x)$, where $\lambda \in \mathbb {Q}$, and the limit of the function $f(x)$ is finite and nonzero at the point $x = x_0$. Petrovi\'c considered the case when $\lambda \neq 0$, i.e. he investigated the question of {\it movable} zeros and poles (generally speaking, critical) of the general solution. He proved the following remarkable statements (see Th\'{e}or\`{e}me  II and Th\'{e}or\`{e}me  V on pages 16-25 in \cite{Petrovich1}).
\smallskip

The first statement: {\it the presence of an inclined edge of the polygon of the equation \eqref {eq1} with the internal normal vector $(\lambda, -1)$, $\lambda <0 $, is a necessary and sufficient condition for the general solution of this equation to have a movable pole of order $-\lambda $}. The second statement can be considered as a consequence of the first: {\it the presence of an inclined edge of the polygon of the equation \eqref {eq1} with the internal normal vector $(\lambda,-1)$, $\lambda > 0$, is a necessary and sufficient condition that the general solution to this equation has a movable zero of order $\lambda $.}
\smallskip

The polygon of the algebraic differential equation \eqref {eq1} corresponding to any of its nonsingular points is defined by Petrovi\'c \cite{Petrovich1} as the {\it convex up} part of the boundary of the convex hull of the set of points $(M_i, N_i) = (p_i + q_i, q_i)$, $i = 1, \ldots, n$, in the plane $OMN$. A distinctive feature of the polygon of a differential equation of precisely the {\it first} order is that, as for the Newton--Puiseux polygon of an algebraic equation,  to each of its points corresponds exactly one monomial of the sum (\ref {eq1}).
\smallskip

We can conclude that Petrovi\'c proved that: {\it if a polygon has an inclined edge with an inner normal vector $(\lambda, -1)$,
$\lambda \ne0$, then there exists a solution $y = (x-x_0)^{\lambda}f(x)$ of the equation \eqref {eq1} such that $f(x_0) \neq 0$, $x_0$ being a parameter.}

\section{The existence of solutions in the form of formal Puiseux series}

Let $x_0$ be a nonsingular point of the equation \eqref {eq1}.
Suppose that the polygon of that equation has an inclined edge $I_{\lambda}$ with an inner normal vector $(\lambda,-1)$, $\lambda \ne0$. Consider the {\it approximate} equation
\begin {equation}\label{ap}
\hat F _ {x_0,\lambda}(y,y ') = \sum_{i:\,Q_i\in I_{\lambda}}a_i (x_0)\,y^{p_i} (y')^{q_i} = 0,
\end {equation}
where only the monomials of the function $F$ corresponding to the points $Q_i=(M_i,N_i)=(p_i+q_i,q_i)$ of this edge participate in the sum.
\medskip

{\bf Definition 2.} A polynomial
$$
P_{x_0,\lambda}(c)=\sum_{i:\,Q_i\in I_{\lambda}}a_i(x_0)\,\lambda^{N_i}c^{M_i}
$$
will be referred to as the {\it characteristic polynomial} corresponding to the edge $I_{\lambda}$.
\medskip

First we note that each nonzero root $c$ of the characteristic polynomial induces a solution $\varphi_0=c(x-x_0)^{\lambda} \ne 0$ of the equation (\ref{ap}). Indeed, if we substitute $\varphi_0$ in the approximate equation \eqref{ap} (in the polynomial $\hat F_ {x_0,\lambda}(y, y')$) we will obtain:
$$
\hat F_{x_0,\lambda}(\varphi_0,\varphi'_0)=P_{x_0,\lambda}(c)\,(x-x_0)^{\gamma},
$$
where $\gamma$ is the (constant) value of the linear function $L_{\lambda}(X, Y)=\lambda X-Y$ on the edge $I_{\lambda}$.
\medskip

{\bf Definition 3.}
We will call the solution $\varphi_0$ {\it non-exceptional} if $c$ is a simple root of $P_{x_0,\lambda}$.
\medskip

Note that the polynomial $P_{x_0,\lambda}$ necessarily has a nonzero root, since it contains at least two monomials and all $M_i$ are different. Hence, the approximate equation (\ref{ap}) necessarily has a nonzero solution of the form $\varphi_0=c(x-x_0)^{\lambda}$.
\medskip

{\bf Theorem 1.} {\it Let $x_0$ be a nonsingular point of the equation \eqref {eq1}. Suppose that the polygon of that equation has an inclined edge $I_{\lambda}$ with an inner normal vector $(\lambda,-1)$, $\lambda=r/s\ne0$.

Each {\rm non-exceptional} solution $\varphi_0=c(x-x_0)^{\lambda}\ne 0$ of the approximate equation \eqref {ap} corresponding to $I_{\lambda}$ is the first term of a uniquely defined Puiseux series $\varphi$ satisfying the equation
\eqref {eq1} and having a nonzero radius of convergence. Moreover the generator of power exponents of $\varphi$ is entirely determined by $\lambda$ and equal to $1/s$.}
\medskip

{\bf Remark 1.} Generally speaking, there can exist several significantly different such solutions of the form $\varphi_0 = c(x-x_0)^{\lambda}\ne 0$ of the equation \eqref {ap}. Here the significant difference means that they are  not transforming into each other after analytic continuation around the point $x = x_0.$ Each such solution will generate its own Puiseux series solution of the original equation.
\medskip

{\bf Remark 2.} The case of $\lambda=0$ ($\varphi_0=c$) is not covered by Theorem 1 and needs an additional study.
\medskip

{\bf Proof}. By technical reason connected with the further application of Malgrange's convergence theorem \cite{Ma},
we will consider the equation (\ref{eq1}) rewritten in the form
\begin{equation}\label{eq1tilde}
\widetilde F(x,y,\delta y)=\sum_{i=1}^na_i(x)(x-x_0)^{m-q_i}\,y^{p_i}(\delta y)^{q_i}=0,
\end{equation}
where $\delta=(x-x_0)(d/dx)$ and $m=\max_i q_i$.

Make the change of variable in the above equation (\ref{eq1tilde}):
$$
y=\varphi_0+(x-x_0)^{\lambda}u, \qquad \delta y=\delta\varphi_0+(x-x_0)^{\lambda}(\delta+\lambda)u.
$$
Then the Taylor formula yields
\begin{eqnarray*}
\widetilde F(x,y,\delta y)&=&\widetilde F(x,\varphi_0,\delta\varphi_0)+(x-x_0)^{\lambda}\,\frac{\partial\widetilde F}{\partial y}(x,\varphi_0,\delta\varphi_0)u+
(x-x_0)^{\lambda}\,\frac{\partial\widetilde F}{\partial(\delta y)}(x,\varphi_0,\delta\varphi_0)(\delta+\lambda)u+ \\
& & +\sum_{k+l\geqslant2}\frac{(x-x_0)^{(k+l)\lambda}}{k!\,l!}\,\frac{\partial^{k+l}\widetilde F}{\partial y^k\,\partial(\delta y)^l} (x,\varphi_0,\delta\varphi_0)u^k((\delta+\lambda)u)^l.
\end{eqnarray*}
Further we note that
\begin{eqnarray*}
\widetilde F(x,\varphi_0,\delta\varphi_0)&=&\sum_{i=1}^na_i(x)\,\lambda^{N_i}c^{M_i}(x-x_0)^{m+\lambda M_i-N_i}=
\Bigl(\sum_{i:\,Q_i\in I_{\lambda}}a_i(x_0)\,\lambda^{N_i}c^{M_i}\Bigr)(x-x_0)^{m+\gamma}+\\
&+&o((x-x_0)^{m+\gamma})=P_{x_0,\lambda}(c)(x-x_0)^{m+\gamma}+o((x-x_0)^{m+\gamma})=
o((x-x_0)^ {m+\gamma}).
\end{eqnarray*}
This is due to the fact that the value $\lambda M_i-N_i=\gamma$ of the function $L_{\lambda}(X,Y)=\lambda X-Y$ is the same for each point $Q_i=(M_i,N_i)$ of the edge $I_{\lambda}$ (recall that $M_i=p_i+q_i$ and $N_i=q_i$). Similarly,
\begin{eqnarray*}
\frac{\partial\widetilde F}{\partial y}(x,\varphi_0,\delta\varphi_0)&=&\sum_{i=1}^na_i(x)\,p_i\,\lambda^{N_i}c^{M_i-1}(x-x_0)^{m+\lambda M_i-N_i-\lambda}=\\
&=&\Bigl(\sum_{i:\,Q_i\in I_{\lambda}}a_i(x_0)\,p_i\,\lambda^{N_i}c^{M_i-1}\Bigr)(x-x_0)^{m+\gamma-\lambda}+o((x-x_0)^{m+\gamma-\lambda}),
\end{eqnarray*}
and since $p_i=(1-\lambda)M_i+\gamma$ for each $Q_i\in I_{\lambda}$, we have
$$
\frac{\partial\widetilde F}{\partial y}(x,\varphi_0,\delta\varphi_0)=(1-\lambda)P'_{x_0,\lambda}(c)(x-x_0)^{m+\gamma-\lambda}+o((x-x_0)^{m+\gamma-\lambda}).
$$
In an analogous way,
\begin{eqnarray*}
\frac{\partial\widetilde F}{\partial(\delta y)}(x,\varphi_0,\delta\varphi_0)&=&\Bigl(\sum_{i:\,Q_i\in I_{\lambda}} a_i(x_0)\,q_i\,\lambda^{N_i-1}c^{M_i-1}\Bigr)(x-x_0)^{m+\gamma-\lambda}+o((x-x_0)^{m+\gamma-\lambda})= \\
&=&P'_{x_0,\lambda}(c)(x-x_0)^{m+\gamma-\lambda}+o((x-x_0)^{m+\gamma-\lambda}),
\end{eqnarray*}
as $q_i=\lambda M_i-\gamma$ for each $Q_i\in I_{\lambda}$. Finally,
$$
(x-x_0)^{(k+l)\lambda}\,\frac{\partial^{k+l}\widetilde F}{\partial y^k\,\partial(\delta y)^l}(x,\varphi_0,\delta\varphi_0)=O((x-x_0)^{m+\gamma}).
$$

Thus after dividing by $(x-x_0)^{m+\gamma}$, the equation (\ref{eq1tilde}) in the unknown variable $u$ takes the form
$$
(1-\lambda)P'_{x_0,\lambda}(c)u+P'_{x_0,\lambda}(c)(\delta+\lambda)u=(x-x_0)^{1/s}\bigl(f(x)+g(x)u+h(x)(\delta+\lambda)u\bigr)+L_2(x,u,(\delta+\lambda)u),
$$
where the positive integer $s$ is the denominator of the exponent $\lambda$: $\lambda=r/s$ ($r$, $s$ being coprime), the coefficients $f,g,h$ of the linear part on the right hand side are polynomials in $(x-x_0)^{1/s}$, the function $L_2$ is polynomial in $(x-x_0)^{1/s}$, $u$, $(\delta+\lambda)u$ and contains only terms at least quadratic in $u$, $(\delta+\lambda)u$. Simplifying the left hand side of the above equation we obtain
$$
P'_{x_0,\lambda}(c)(\delta+1)u=(x-x_0)^{1/s}\bigl(f(x)+g(x)u+h(x)(\delta+\lambda)u\bigr)+L_2(x,u,(\delta+\lambda)u).
$$
Making the change of the independent variable, $x-x_0=t^s$, we come to the analytic ODE near $0\in{\mathbb C}^3$:
\begin{equation}\label{eq3}
P'_{x_0,\lambda}(c)((1/s)\delta_t+1)u=t\bigl(\tilde f(t)+\tilde g(t)u+\tilde h(t)((1/s)\delta_t+\lambda)u\bigr)+\widetilde L_2(t,u,((1/s)\delta_t+\lambda)u),
\end{equation}
where $\delta_t=t(d/dt)$.

The equation (\ref{eq3}) has a unique power series solution $\hat u=\sum_{k\geqslant1}c_k\,t^k\in{\mathbb C}[[t]]$, where
$$
P'_{x_0,\lambda}(c)(1+1/s)c_1=\tilde f(0),
$$
and other $c_k$'s with $k>1$ are uniquely determined by the previous ones $c_1,\ldots,c_{k-1}$. By the Malgrange theorem \cite{Ma}, this series has a nonzero radius of convergence. Hence the initial equation (\ref{eq1}) possesses the unique Puiseux series solution $\varphi=c(x-x_0)^{r/s}+\sum_{k\geqslant1}c_k\,(x-x_0)^{(r+k)/s}$ starting with the term $\varphi_0=c(x-x_0)^{r/s}$ and having a nonzero radius of convergence. This proves Theorem 1. {\hfill $\Box$}

\section{From formal to actual solutions}

Now Theorem 1 allows us to prove the statement on the convergence of formal Puiseux series solutions of a general algebraic ODE of first order.
\medskip

{\bf Theorem 2.} {\it All formal Puiseux series in powers of $x-x_0$ satisfying the equation \eqref {eq1} and starting with a {\rm non-exceptiona}l term $c(x-x_0)^{\lambda}$, $\lambda\ne0$, where $x_0$ is a nonsingular point of the equation, converge in some small neighborhood of it.}
\medskip

{\bf Proof.} From the initial arguments in the proof of Theorem 1, carried out in the reverse order, it follows that the first term $\varphi_0 = c(x-x_0)^{\lambda}$, $\lambda\in{\mathbb Q}^*$, of the formal Puiseux series $\varphi$ satisfying the equation (\ref {eq1}), is a solution of the approximate equation \eqref {ap}, corresponding to the inclined edge of the polygon of the equation (\ref {eq1}) with the inner normal vector $(\lambda, -1)$. By virtue of uniqueness and convergence of the Puiseux series starting with $\varphi_0$ and satisfying the equation (\ref {eq1}), the series $\varphi$ has a nonzero radius of convergence. {\hfill $\Box$}
\medskip

Let us conclude with several examples illustrating that, although the Painlev\'e theorem and Theorems 1, 2 have an intersection in applications they naturally complement each other on the important questions of convergence.
\medskip

{\bf Example 1.} Let us consider the first example, an ODE for which the Painlev\'e--Picard system \eqref{systP} is compatible for any $x\in\mathbb C$:
\begin{equation}\label{exp2b1}
F(x,y,y')=y'^3-(y-x^4)^2=0.
\end{equation}
It follows that the fixed finite set $\Sigma$ where solutions could have singular points $x=x_0$ of non-algebraic type and thus formal Puiseux series solutions in powers of $x-x_0$ could be divergent, is somehow hidden inside the equation. However, one can easily see that for any point $x_0\ne0$ the polygon of the equation consists of one edge connecting the points $(0,0)$ and $(3,3)$ thus having the inner normal vector $(1,-1)$. The characteristic polynomial corresponding to this edge $I_{\lambda=1}$ is
$$
P_{x_0,1}(c)=c^3-x_0^8,
$$
and all its roots are simple ($x_0\ne0$). Hence Theorem 1 gives that for each of the three values of $x_0^{8/3}$, there is a unique formal Puiseux series solution of \eqref {exp2b1} beginning with $x_0^{8/3}(x-x_0)$ and having a nonzero radius of convergence (this is a Taylor series in fact).

One could expect the singular point $x=0$ of the equation to be a candidate for a non-algebraic singular point of the general solution of (\ref{exp2b1}). Actually at this point there is a formal Puiseux series solution $y=x^4\sum_{k\geqslant0}a_kx^{k/2}=x^4+8x^{9/2}+108x^5+1863x^{11/2}+37665x^6+\dots$. The numerical evidences indicate  that this series is rapidly diverging.
\medskip

{\bf Example 2.} The next example shows that in the case where the set $\Sigma$ is simply detected by an equation, the convergence of some Puiseux series solutions in the variable $x-x_0$, $x_0\in\Sigma$, can still be established by Theorems 1, 2.

Consider the equation
\begin{equation}\label{exeq1}
F(x,y,y')=-2y^3y'^3-y^2y'^2+y-(x^3+1)=0.
\end{equation}
In this case the Painlev\'e--Picard system (\ref{systP}) is compatible for $x=0$ and this point thus being an element of $\Sigma$ could be the fixed singular point of non-algebraic type for the general solution. The equation possesses a formal power series solution
\begin{equation}\label{ex1}
1+\sum_{k\geqslant3}c_k(k-2)!\,x^k
\end{equation}
in powers of $x$, starting with the constant term $c_0=1$ and having a zero radius of convergence, which says
that there indeed exist solutions with the singularity of non-algebraic type at $x=0$. Indeed, after the change $y=1+w$ of the dependent variable, we come to the equation
$$
(1+w)^2(1+2(1+w)w')w'^2=w-x^3.
$$
By substituting into the above equation $w=\sum_{k\geqslant3}c_k(k-2)!\,x^k$, we get $c_3=1$ and, further,
$$
(1+6x^2+\ldots)w'^2=\sum_{k\geqslant4}c_k(k-2)!\,x^k,
$$
hence
$$
(1+6x^2+\ldots)\Bigl(\sum_{l\geqslant3}c_l(l-2)!\,l\,x^{l-1}\Bigr)\Bigl(\sum_{m\geqslant3}c_m(m-2)!\,m\,x^{m-1}\Bigr)=
\sum_{k\geqslant4}c_k(k-2)!\,x^k.
$$
Thus, $c_4=9/2$. By comparing the coefficient $c_k(k-2)!$ of $x^k$ on the right hand side with the coefficient of $x^k$ from the left hand side, we get
$$
c_k(k-2)!=6c_{k-1}(k-3)!(k-1)+P_k(c_3,\ldots,c_{k-2}), \qquad k=5, 6, \ldots,
$$
where $P_k$ is a polynomial with {\it positive} coefficients. From the relation
$$
c_k=\frac{6(k-1)}{k-2}\,c_{k-1}+\frac1{(k-2)!}\,P_k(c_3,\ldots,c_{k-2}), \qquad k=5, 6, \ldots,
$$
the coefficients $c_k$ are uniquely determined. Taking into account that  $c_3=1$ and $c_4=9/2$, by using the mathematical induction, we prove $c_k\geqslant 1$ for all $\;k\geqslant3$. Consequently, the series (\ref{ex1})
has the radius of convergence equal to zero.

On the other hand, the polygon of the equation (\ref{exeq1}) corresponding to its nonsingular point $x=0$ consists of one edge connecting the points $(0,0)$ and $(6,3)$ and also containing the point $(4,2)$. The inner normal vector of this edge is
$(1/2,-1)$ and the corresponding characteristic polynomial is
$$
P_{0,1/2}(c)=-\frac14\,c^6-\frac14\,c^4-1,
$$
which has six simple roots $\pm c_i$, $i=1,2,3$. Therefore we have three essentially different convergent Puiseux series solutions of (\ref{exeq1}) in the variable $x$, beginning with $c_ix^{1/2}$.
\medskip

Concerning formal Puiseux series solutions in the variable $x-x_0$, where $x_0$ is a nonsingular point of an equation, that
begin with an {\it exceptional} term $c(x-x_0)^{\lambda}$, $\lambda\ne0$, they can also be convergent.
\medskip

{\bf Example 3.} Consider the equation
$$
(y'-1)^2-9x=0.
$$
For its nonsingular point $x_0=0$ and the edge $I_{\lambda=1}$ of the polygon, one has an exceptional solution $\varphi_0=x$ of the corresponding approximate equation $(y'-1)^2=0$ (since $c=1$ is a double root of the corresponding
characteristic polynomial $P_{0,1}(c)=(c-1)^2$). However, this solution is extended to an actual solution $x+2x^{3/2}$ of the initial equation. Though, note that the generator $1/2$ of the power exponents of the last truncated Puiseux series is not
determined by that of the starting monomial $x$.
\medskip

But in general the condition of non-exceptionality is essential and cannot be omitted, as we see in the following example.
\medskip

{\bf Example 4.} Let us consider the differential equation
$$
F(x,y,y')=y'^2+2y'-y+ (1-x+x^3)=0.
$$
The point $x_0=0$ is a nonsingular point of this equation, the polygon has an inclined edge $I_{\lambda=1}$ containing the points $(0,0)$, $(1,1)$ and $(2,2)$. The corresponding characteristic polynomial is $P_{0,1}(c)=(c+1)^2$. The equation has
a formal power series solution
$$
-x+x^3+9x^4+216x^5+7776x^6+\ldots,
$$
beginning with the exceptional term $\varphi_0=-x$ (since $c=-1$ is a double root of $P_{0,1}(c)$). This formal power series has a zero radius of convergence indeed, which can be proved by making the change $y=-x+x^3+w$ of the dependent variable and applying reasoning similar to that from Example 2, to the transformed equation
$$
w'^2+6x^2w'-w+9x^4=0.
$$

\subsection*{Acknowledgements}
We thank the referee  of the previous journal version of this paper for providing us with the last example, which helped us refine and finalize this version of the paper.
This research has been partially supported by the Simons Foundation grant no. 854861.

\end{document}